\documentclass[preprint,12pt]{elsarticle}

\usepackage{amssymb, amsmath}
\usepackage{soul}

\def\RR{{\mathbb R}}
\def\NN{\mathbb N}
\def\ZZ{\mathbb Z}

\newtheorem {theorem} {Theorem}
\newtheorem {lemma} {Lemma}

\newtheorem {definition} {Definition}

\begin{document}

\begin{frontmatter}


 \author{{Vladyslav Babenko}}
 \ead{babenko.vladislav@gmail.com}
 \ead[url]{http://math.kennesaw.edu/~ybabenko/VF$\_$Babenko}
 \address{Department of Mathematical Analysis and Theory of Functions, Dnepropetrovsk National University, pr. Gagarina 72, Dnepropetrovsk, Ukraine 49050}
 
 \author{Yuliya Babenko\corref{cor1}}
 \ead{ybabenko@kennesaw.edu}
 \ead[url]{http://math.kennesaw.edu/~ybabenko/}
 \cortext[cor1]{Corresponding author}
 \address{{Department of Mathematics and Statistics,
Kennesaw State University,
1000 Chastain Road, 1601,
Kennesaw, GA, USA 30144 }}
 
 \author{{Oleg Kovalenko}}
 \ead{olegkovalenko90@gmail.com}
 \address{Department of Mathematical Analysis and Theory of Functions, Dnepropetrovsk National University, pr. Gagarina 72, Dnepropetrovsk, Ukraine 49050}

\title{Kolmogorov's Problem on the Class of Multiply Monotone Functions}




\begin{abstract}
{ In this paper we give necessary and sufficient conditions
for the system of positive numbers $ M_{k_1}, M_{k_2},...,
 M_{k_{d}},$ $0\leq k_1<...<k_{d} {\leq} r$, to guarantee the
existence of an $r$-monotone function defined on the negative
half-line $\RR_-$ and such that $ \|x^{(k_i)}\|_{\infty}=M_{k_i},\;
i=1,2,...,d. $}
\end{abstract}


\end{frontmatter}
\section {Introduction} \label{Intro}

 Let the domain $G$ be the real line $\RR=(-\infty , \infty )$ or the negative half-line $\RR_-=(-\infty , 0]$. In addition, let $L_{\infty}(G)$ be the space of measurable essentially
bounded functions $x:G\to \RR$ with usual norm $\| \cdot \|=\|
\cdot \|_{L_{\infty}(G)}$.

 For $r\in \NN$ we denote by $L^r_{\infty} (G)$ the space of functions $x:G\to \RR$ that have
locally absolutely continuous derivative of order $r-1$, $x^{(0)}=x$
(in the case $G=\RR_-$
we take, as usual, the one-sided derivative at the point $x=0$),
and such that $x ^ {(r)} \in L_{\infty} (G)$. Let
$L^r_{\infty,\infty}(G)=L^r_{\infty}(G)\cap L_{\infty}(G)$.

Kolmogorov (see   ~\cite{Kol1},~\cite{Kol2},~\cite{Kol3})
formulated
the following problem:

{\bf Kolmogorov's problem.}
{\it Let some class of functions $X\subset L^r_{\infty,\infty}(G)$ and
 an arbitrary system of $d$ integers
$$
0\leq k_1 < k_2 < ... < k_{d}{\leq}r
$$
be given. The problem is to find the necessary and sufficient conditions on the
system of positive numbers
\begin{equation}\label{ms}
M_{k_1}, M_{k_2}, ..., M_{k_{d}}
\end{equation}
to guarantee the existence of a function $x \in X$
 such that
\begin{equation}\label{co}
\|x^{(k_i)}\|=M_{k_i}, \;\;\;\;\; i=1,...,d.
\end{equation}
}
 In ~\cite{Kol1},~\cite{Kol2} (see also ~\cite{Kol3}) Kolmogorov stated and solved this
 problem for the case $d=3$, $X= L^r_{\infty,\infty}(\RR)$, {$k_1=0$ and $k_d=r$}. Rodov was the first who considered
 Kolmogorov's problem for $d>3$ (see ~\cite{Rod}).
 Note that for $d=2$ and $X\subset L^r_{\infty,\infty}(G)$
 (invariant with respect to multiplication of a function and its argument by an
 arbitrary positive constant) the solution to Kolmogorov's problem is trivial:
 for any pair of positive numbers $M_{k_1}$ and $M_{k_2}$ there exists a function
 $x\in X$ such that (\ref{co}) are satisfied.


We point out that for $d\geq 3$, any sharp inequality for the norms of consecutive derivatives $\|x^{(k_i)}\|$, $\|x^{(k_l)}\|$ and $\|x^{(k_m)}\|$ ($k_i<k_l<k_m$) of the function $x\in X$ provides a necessary condition that numbers $M_{k_i}$, $M_{k_l}$, and $M_{k_m}$ from (\ref{ms}) must satisfy in order to guarantee the existence of function $x\in X$ such that all conditions (\ref{co}) are satisfied. Such inequalities are called Landau-Kolmogorov type inequalities. In turn,
Landau-Kolmogorov type inequalities have numerous applications in Approximation Theory, in particular in the approximation of unbounded operators by bounded ones (Stechkin's problem), embedding theorems, Bernstein type inequalities, Nikol'skii type inequalities for trigonometric polynomials and splines (with different metrics), optimization of quadrature formulae, and investigating extremal properties of polynomials and splines. For more details on existing results, applications, open problems and further references we direct interested readers to ~\cite{Are96}, Chapters 7 and 8 in ~\cite{BKKP},  ~\cite{B}, ~\cite{BKP}, ~\cite{BKP2}. Note that, in the case $d=3$, such sharp inequalities  usually provide a complete solution (necessary and sufficient conditions) to Kolmogorov's problem.

Clearly, Kolmogorov's problem can be formulated for various classes of functions $X$, for functions
defined on other domains (finite interval, circle, etc.), for multivariate functions, for other norms (potentially different
for different derivatives). In particular, Kolmogorov's problem has been considered for the class of multiply monotone functions (see ~\cite{us, ol}) which is the class of interest in this paper (and defined later), and the class of absolutely monotone functions ~\cite{us_abs}.


The paper is organized as follows. Section \ref{known} presents the
cases of the known solutions to Kolmogorov's problem and discusses what could naturally be considered a solution to
Kolmogorov's problem. Section \ref{main} contains necessary
definitions and statements of the main results, which provide the
answer to Kolmogorov's problem for arbitrary  $d$ positive
numbers to be values of norms of prescribed derivatives of an
$r$-monotone function from $L^r_{\infty,\infty}(\RR_-)$.
Section \ref{auxiliary} is dedicated to the proofs of several supporting results which, when combined,
constitute the proof of the main result of the paper.

 \section{Review of the known results}\label{known}


The complete solution for three numbers (more precisely, in the case $d=3$, $k_1=0$, $k_3=r$ and $L_{\infty,\infty}^r(\RR)$)
%
was given by Kolmogorov (in the case $r=2$ and $k=1$ it follows from the results of ~\cite{had}, for
all cases with $r<5$ and for $r=5$ and $k_2=2$  it follows from ~\cite{bosse}) , who showed that for any three positive numbers $M_0, M_k, M_r$, $k_1=0<k<r$, there exists a function $x \in
L_{\infty,\infty}^r(\RR)$ that has these numbers as values of the
norm of the function itself, its $k^{\rm{th}}$ and its $r^{\rm{th}}$
derivatives, respectively, if and only if
$$
M_k\leq
\frac{\left\| \varphi _{r-k}\right\| }{%
\left\| \varphi _r\right\| ^{1-k/r}}%
M_0^{1-k/r}M_r^{k/r},
$$
or, that is equivalent,
$$
M_0\geq
\frac{\|\varphi_r\|}{\|\varphi_{r-k}\|^{\frac{r}{r-k}}}M_k^{\frac{r}{r-k}}M_r^{-\frac{k}{r-k}},
$$
 where
 $\varphi _r$ is the $r^{\rm{th}}$ periodic integral with zero mean value on a period
of the function $\varphi_0\left( t\right) ={\rm sgn}(\sin t)$. Note
that the case when $k_1>0$ was resolved in ~\cite[Chapter 9]{BKKP}.


  The only other known solutions for the case when the domain is the whole real
  line $\RR$ are in the following cases:
\begin{enumerate}
\item $X=L^{r}_{\infty,\infty}({\RR})$; $k_1=0,\; k_2=r-2,\; k_3=r-1,\; k_4=r$ (Rodov 1946~\cite{Rod}).
\item $X=L^{r}_{\infty,\infty}({\RR})$; $k_1=0<k_2<k_3=r-2,\; k_4=r-1, \;k_5=r$ (Rodov {1946}~\cite{Rod}).
\item $X=L^{r}_{\infty,\infty}({\RR})$; $k_1=0<k_2< k_3=r-1,\; k_4=r$ (Dzyadyk, Dubovik 1975~\cite{DD2}).
\item $X=L^{r}_{\infty,\infty}({\RR})$; $k_1=0<k_2< k_3=r-2,\;k_4=r$
(V. Babenko, Kovalenko 2012~\cite{BK})
\end{enumerate}
In \cite{Rod2} Rodov found sufficient conditions for {$r=5$ and} sets
$$M_0, M_1, M_2, M_5\;\; {\rm and}\;\; M_0, M_1, M_2,M_3, M_4, M_5.$$

For the general case of arbitrary $d$ {(with $k_d=r$)}, the only result is by Dzyadyk and Dubovik~\cite{DD1} and it
provides a sufficient condition for the function $x \in
L_{\infty,\infty}^r(\RR)$ to exist.

In the case when the domain is the half-line $\RR_-$, solutions to
Kolmogorov's problem are known in the following cases:
\begin{enumerate}
\item $X=L^{r}_{\infty,\infty}({\RR}_-)$, $k_1=0<k_2< k_3=r$ (partial cases follow from
the results of Landau 1913~\cite{Lan} and Matorin 1955~\cite{Mat}.
The complete solution, although in implicit form, follows from the
work of Schoenberg and Cavaretta 1970,~\cite{ShC1, ShC2}).
\item $X=L^{r}_{\infty,\infty}({\RR}_-)$, $k_1=0<k_2<k_3=r-1,\;k_4=r$ (V.~Babenko and Britvin 2002~\cite{BB}).
\end{enumerate}

Given $r,m\in \mathbb{Z}_+, \;m\le r$, denote by
$L_{\infty,\infty}^{r,m}(\RR_-)$ the class of functions $x \in
L^{r}_{\infty,\infty}(\RR_-)$ that are nonnegative along with all
their derivatives up to and including order $m$ (derivative of order
$m$ must be nonnegative almost everywhere if $m=r$). We will henceforce call
it the class of $m$-monotone functions.

In 1951, Olovyanishnikov~\cite{ol} obtained {an
elegant} solution to Kolmogorov's problem for three numbers and
for the class of multiply monotone functions defined on a half-line.
He showed that for arbitrary $k,r\in \NN $, $k<r$, and for any three
positive numbers $M_0, M_{k}, M_r$, there exists a function $x\in
L^{r,r-1}_{\infty,\infty}(\RR_-)$ such that
$$
\|x\|=M_0,\;\;\|x^{(k)}\|=M_{k},\;\;\|x^{(r)}\|=M_{r},
$$
if and only if these numbers satisfy Kolmogorov type inequality
$$
M_0\geq
\frac{\|\phi_r\|}{\|\phi_{r-k}\|^{\frac{r}{r-k}}}M_k^{\frac{r}{r-k}}M_r^{-\frac{k}{r-k}},
$$
where $\phi_r(t):=\displaystyle \frac{l(t+a)_+^r}{r!}$, which in
particular provides an elegant sharp constant in the above inequality.

In
~\cite{subbchern} and independently in ~\cite{usEJA}, this result was extended to $(r-2)$ -monotone functions (moreover, the last one was in the case of arbitrary norms).
In addition, in ~\cite{us} the solution was obtained for the class $L^{r,r}_{\infty,\infty}({\RR_-})$  in the case $k_1>0$.

Other known solutions for the classes of multiply monotone functions
defined on a half-line $\RR_-$ are in the following cases:
\begin{enumerate}
\item $X=L^{r,r-2}_{\infty,\infty}({\RR}_-)$ and $k_1=0< k_2< k_3=r-1,\; k_4=r$ (Yattselev 1999~\cite{Yat}). Note that the order of only one derivative ($k_2$) is not fixed.
\item $X=L^{r,r}_{\infty,\infty}({\RR}_-)$ and $k_1=0< k_2< k_3<k_4=r$ (V. Babenko, Y. Babenko 2007~\cite{us}). Note that in this case the orders of  two intermediate derivatives ($k_2$ and $k_3$) are arbitrary.
\end{enumerate}

As mentioned above, in this paper we give the complete solution (for an arbitrary number
$d$ of prescribed derivatives) to Kolmogorov's problem for the class $L_{\infty ,
\infty}^{r,r}(\mathbb{R}_-)$ of $r$-monotone functions.



Note that in all the previously described cases the solution can have an alternative interpretation.  In order to introduce it, we need the following notation. Let $d\in\NN$ and integers $0\leq k_1<k_2<...<k_{d}{\leq}r$ be given. Let
${\bf k}:=(k_1,\ldots ,k_{d})$. For $i=1,\dots ,d$ set
${\bf{k}}\,^i=(k_i,k_{i+1},\ldots ,k_{d})$, so that ${\bf
k}\,^1={\bf k}$. The sets of positive numbers $\{M_{k_1},\ldots
,M_{k_{d}}\}$ and  $\{M_{k_i},\ldots ,M_{k_{d}}\}$ we will denote by
$M_{{\bf k}}$ and $M_{{\bf k}^i}$, respectively. For a given vector ${\bf k}:=(k_1,\ldots ,k_{d})$ and function $x\in X$ we set
$$
M_{\bf k}(x):=(M_{k_1}(x),\ldots ,M_{k_d}(x)),
$$
where
$$
M_{k_i}(x)=\| x^{(k_i)}\|,\;\; i=1,\ldots ,d.
$$

\begin{definition}\label{admissible}
We say that the set $M_{{\bf k}}$ is {\bf admissible} for the class
$X\subset L^r_{\infty,\infty}(G)$, if there exists a function $x\in
X$ such that $\left\|x^{(k_i)}\right\|=M_{k_i}$, $i=1,2,...,d$ (or,
for short, $M_{{\bf k}}(x)=M_{{\bf k}}$). We denote by
{$A_{d}(X)=A_d(X,{\bf k})$} the collection of all
admissible sets $M_{{\bf k}}$. Notation $M_{{\bf k}}\in A_{d}(X)$
{{means}} that the set $M_{{\bf k}}$ is admissible
for the class $X$.
\end{definition}

So in the above cases, solution to Kolmogorov's problem can be interpreted as follows. For a wide class of functions
 $X$ one finds a $d$-parametric family $F$ of functions, such that
 \begin{equation}\label{pa}
 A_{d}(X)=\{ M_{{\bf k}}(x)\; :\; x\in F\}.
 \end{equation}
If we interpret the set $A_{d}(X)$ in such a way, it is natural to assume the minimality requirement on the set: for example, to require that the removal of at least one function from  the set $F$ would not generate the whole set  $A_{d}(X)$.

\begin{definition}
Minimal $d$-parametric family of functions $F\subset X$ such that (\ref{pa}) is true is called the {\bf generating set} for $A_d(X,{\bf k})$. We denote it by $F_d(X,{\bf k})$.
\end{definition}

With the help of this definition we can re-formulate Kolmogorov's problem.

\bigskip

{\bf Kolmogorov's problem.}
{\it For a given class of functions $X\subset L^r_{\infty,\infty}(G)$, fixed $d$ and {vector
${\bf k}$}, find (or characterize) the set $F_{d}(X,{\bf k})$ for $A_d(X,{\bf k})$.}

\bigskip

From this perspective, the solution to classical Kolmogorov's problem (for three numbers) may be interpreted as follows. Given numbers
$M_k$ and $M_r$, we can find values of parameters
 $\lambda >0$, $\alpha>0$, and set $\varphi
_{\lambda, r}(t):=\lambda ^{-r}\varphi_r(\lambda t)$ such that we obtain
 $$
 M_k(\alpha\varphi _{\lambda,r})=M_k, \qquad  M_r(\alpha\varphi _{\lambda,r})=M_r.
 $$
Then the desired function is
$$
\alpha\varphi _{\lambda,
r}+C,
$$
where $ \alpha>0,\lambda >0$ and $C=M_0-\frac{\|\varphi_r\|}{\|\varphi_{r-k}\|^{\frac{r}{r-k}}}M_k^{\frac{r}{r-k}}M_r^{-\frac{k}{r-k}}\geq 0$.
Therefore,
$$
F_3(L_{\infty,\infty}^r({\RR}), {\bf k})=\{ \alpha\varphi _{\lambda,
r}+C\; :\; \alpha>0,\lambda >0, C\ge 0\},
$$ where
${\bf k}=(k_1=0, k_2=k, k_3=r).
$

In the view of this framework, Olovyanishnikov's results mentioned above (see ~\cite{ol})
can be re-written as
$$
F_3(L^{r,r-1}_{\infty,\infty}({\RR_-}), {\bf k})=F_3(L^{r,r}_{\infty,\infty}({\RR_-}), {\bf k})=\{\alpha\phi_r(\lambda t) +C: \alpha>0,\;\, \lambda>0, \;\;{C\geq 0}\},
$$
where $ {\bf k}=(k_1=0, k_2=k, k_3=r).$

\section{Statements of the main results} \label{main}

In order to state the main results of this paper we need the following
definitions.

Let $r,s\in \NN$, $a_1>a_2>...>a_s>0$ and $l>0$ be given. We
define a function, later referred to as a spline of order $r$ with
knots { $-a_1<-a_2<...<-a_s<0$}, as follows
\begin{equation}\label{phi}
\phi(a_1,a_2,...,a_s,l;t):=\frac l {r!}\sum\limits_{j=1}^s
(-1)^{j+1}(t+a_j)^r_+.
\end{equation}

%

Note that, since this function and all its derivatives up to order $r-1$
are monotone on $\RR_-$, the uniform norms of all derivatives of
order less than $r$ are achieved at zero and, hence, are equal to
\begin{equation}\label{norms}
\|\phi^{(k)}(a_1,a_2,...,a_s,l)\|=\frac{l}{(r-k)!}\sum\limits_{j=1}^s
(-1)^{j+1} a_j^{r-k}, \qquad k=0,\dots,r-1,
\end{equation}
$$\|\phi^{(r)}(a_1,a_2,...,a_s,l)\|=l.$$
In addition, let
\begin{equation} \label{Phi}
\Phi_{r,n}:=\left\{\phi(a_1,a_2,...,a_s,l;t):s\in {\NN}, s\leq n,\;
a_1>a_2>...>a_s>0,\; l>0\right\}
\end{equation}
denote the collection of all splines of order $r$ with no more than $n$ knots.

%
%

Let $d\in\NN$ and integers $0\leq k_1<k_2<...<k_{d} {\leq} r$ be given.



As previously mentioned, the goal is to find an admissible set for the set of $r$-monotone functions, i.e. the set $A_d(L^{r,r}_{\infty,\infty}({\RR}_-))$ which we will refer to as $A_d$ for brevity, and its generating set. In addition, we need the following definitions.

\begin{definition}\label{veryType 2}
An admissible set $M_{{\bf k}}\in A_d$  is called of {\bf Type 1}, if
there exists a spline $\phi\in\Phi_{r,d-1}\setminus \Phi_{r,d-2}$
such that
$\left\|\phi^{(k_i)}\right\|=M_{k_i}$, $i=1,2,...,d$.  We denote by $A^1_{d}$ the collection of all admissible sets $M_{{\bf k}}$ of Type 1. 
\end{definition}


\begin{definition}\label{Type 2}
An admissible set $M_{{\bf k}}\in A_d$  is called of {\bf Type 2}, if
there exists a spline $\phi\in\Phi_{r,d-2}$ such that
$\left\|\phi^{(k_i)}\right\|=M_{k_i}$, $i=1,2,...,d$. We denote by $A^2_{d}$ the collection of all admissible sets $M_{{\bf k}}$ of Type 2. 
\end{definition}

%

\begin{definition}\label{bad}
An admissible set $M_{{\bf k}}\in A_d$ with $k_1=0$  is called of {\bf
Type 3}, if there exists a constant $C>0$  and a spline
$\phi\in\Phi_{r, d-1}$ such that
$\left\|(\phi+C)^{(k_i)}\right\|=M_{k_i}$, $i=1,2,...,d$, and this
set
is not a Type 1 set. We denote by $A^3_{d}$ the collection of all admissible sets $M_{{\bf k}}$ of Type 3. 
\end{definition}


%

{The following theorem can be viewed as a generalization of Olovyanishnikov's inequality to $d$ norms of consecutive derivatives of a function.}

\begin{theorem}\label{th1} (Existence {and extremal properties} of a
spline). Let $r,d\in\NN$, $d\geq {3}$,
 and integers $0\leq k_1<...<k_{d} {\leq} r$ be given. Let also $x(t)\in L_{\infty,\infty}^{r,r}(\RR_-)$ be given. Then there exists a spline $\phi(t)=\phi(a_1,a_2,...,a_s,l;t)\in \Phi_{r, d-1}$ and a  number $C\geq 0$ such that
 \begin{equation}
  \left\|x^{(k_i)}\right\|=\left\|(\phi+C)^{(k_i)}\right\|, \qquad i=1,2,...,d.
  \end{equation}
  {   Moreover, if $k_d=r$ and $k\in \ZZ_+$ are such that for some $i = 0,1,\dots,d-1$  we have $k_i<k<k_{i+1}$,  $(k_0:=-1)$ and $x^{(k_1)}\neq \phi^{(k_1)}$, then $$(-1)^i\left\|x^{(k)}\right\|>(-1)^i\left\|\phi^{(k)}\right\|,$$} { and in the case when $k_d<r$ and $x^{(k_1)}\neq \phi^{(k_1)}$ we have $$\left\|x^{(r)}\right\|>\left\|\phi^{(r)}\right\|.$$}
\end{theorem}
The corresponding spline will be denoted by $\phi({\bf k}, x;t)$.
Given $M_{{\bf k}}$, the spline $\phi$ such that
$\|\phi^{(k_i)}\|=\phi^{(k_i)}(0)=M_{k_i},\; i=1,\ldots ,d$, we
shall denote by $\phi (M_{{\bf k}};t)$.

\begin{theorem} (Solution to Kolmogorov's problem in alternative form).\label{T3}
Let $d\in\NN$, $d> 3$, and ${\bf k}:=(k_1,\ldots ,k_{d})$ with
$0\leq k_1<k_2<...<k_{d}{\leq}r$ be given. {Then in the case when
$k_d=r$}
$$
F_d(L_{\infty,\infty}^{r,r}(\RR_-), {\bf k})=\{  \Phi_{r, d-1}  \setminus \Phi_{r,d-2}\}
\bigcup \{ \phi
+C\; :\; \phi \in \Phi_{r, d-2},\,C\geq 0\},
$$
{and in the case when $k_d<r$ $$
F_d(L_{\infty,\infty}^{r,r}(\RR_-), {\bf k})= \Phi_{r, d-1} \setminus \Phi_{r,d-2}.
$$}
\end{theorem}


 {Another way of stating the solution to Kolmogorov's problem in the case when $k_d=r$ is given in the following theorem.}

\begin{theorem}\label{th3} (Solution to Kolmogorov's problem for $k_d=r$)
Let $d\in\NN$, $d{\geq}3$, and $0\leq k_1<k_2<...<k_{d}=r$ be
nonnegative integers. Then
$$
\{ M_{{\bf k}} \in A_{d}(L^{r,r}_{\infty,\infty}({\RR_-}))\}
$$
$$\Longleftrightarrow$$
$$
\left\{\begin{array}{rcl}
M_{{\bf k}^2} \in A^1_{d-1}\;\;\;\; \;\; \; \\
M_{k_1}\ge \|\phi^{(k_1)}(M_{{\bf k}^2})\|\\
\end{array}\right\}\;\bigvee \left\{\begin{array}{rcl}
M_{{\bf k}^2} \in A^2_{d-1} \;\;\;\; \;\; \; \\
k_1>0\;\;\;\; \;\; \;\;\;\;\\
M_{k_1}=\| \phi^{(k_1)}(M_{{\bf k}^2})\| \\
\end{array}\right\}\;\bigvee
\left\{\begin{array}{rcl}
M_{{\bf k}^2} \in A^2_{d-1}\;\;\;\; \;\; \; \\
k_1=0\;\;\;\; \;\; \;\;\;\;\\\
M_{k_1}\ge \|\phi^{(k_1)}(M_{{\bf k}^2})\|\\
\end{array}\right\}.
$$
Moreover,
$$
\left\{\begin{array}{rcl}
M_{{\bf k}^2}\in A^1_{d-1} \;\;\;\;\;\\
M_{k_1}\ge \|\phi^{(k_1)}(M_{{\bf k}^2})\|\\
\end{array}\right\}\Longrightarrow \left\{\begin{array}{rcl}
M_{{\bf k}}\in A^1_{d}\;\;\;\;{\rm if} \;\;\;\; M_{k_1}> \|\phi^{(k_1)}(M_{{\bf k}^2})\| \\
M_{{\bf k}}\in A^2_{d}\;\;\;\;{\rm if} \;\; \;\; M_{k_1}=
\|\phi^{(k_1)}({M_{{\bf k}^2}})\|
\end{array}\right\},
$$
$$
\left\{\begin{array}{rcl}
M_{{\bf k}^2}\in A^2_{d-1} \;\;\;\;\;\ \\
k_1>0\;\;\;\;\;\;\;\;\\
M_{k_1}= \|\phi^{(k_1)}(M_{{\bf k}^2})\| \\
\end{array}\right\}\Longrightarrow \{M_{{\bf k}}\in A^2_{d}\},
$$

$$
\left\{\begin{array}{rcl}
M_{{\bf k}^2}\in A^2_{d-1} \;\;\;\;\;\\\
k_1=0 \;\;\;\;\;\;\;\;\\
M_{k_1}\ge \|\phi^{(k_1)}(M_{{\bf k}^2})\| \\
\end{array}\right\}\Longrightarrow \left\{\begin{array}{rcl}
M_{{\bf k}}\in A^2_{d}\;\;\;\;{\rm if} \;\;\;\; M_{k_1}= \|\phi^{(k_1)}(M_{{\bf k}^2})\| \\
M_{{\bf k}}\in A^3_{d}\;\;\;\;{\rm if} \;\; \;\; M_{k_1}>
\|\phi^{(k_1)}({M_{{\bf k}^2}})\|
\end{array}\right\}.
$$
\end{theorem}

{{\bf{Remark 1.}}}
{It is easy to see that, in the case $d=2$, we have $\left\{M_{k_1}, M_{k_2}\right\}\in A_1^1$ for all $0\leq k_1<k_2\leq r$ and all positive $M_{k_1}, M_{k_2}.$ Hence, Theorem~\ref{th3} in the case $d=3$ can be re-written in the following way: the set
$\left\{M_{k_1}, M_{k_2}, M_r\right\}$ is admissible
if and only if}

\begin{equation}\label{olov}
M_{k_1}\geq \frac{(r-k_2)!^{\frac{r-k_1}{r-k_2}}}{(r-k_1)!}M_{k_2}^{\frac{r-k_1}{r-k_2}}M_r^{\frac{k_1-k_2}{r-k_2}}.
\end{equation}
In the case $k_1=0$ this result is contained in the work of Olovyanishnikov~\cite{ol}. Cases $d=3$,
$k_1>0$ and $d=4$, $k_1= 0$ are contained in the work of V. Babenko
and Y. Babenko ~\cite{us}.

{\bf Remark 2.} Note that Theorem 3 gives, in particular, an ``algorithm'' to
test an arbitrary set of positive numbers $\{M_{k_1},\dots, M_{k_d}\}$
for existence of the needed function from the class $L^{r;r}_{\infty,\infty}(R)$ and, moreover, a way to construct such
a function (when possible). Without technical details the idea is the following. Given the set $\{M_{k_1},\dots, M_{k_d}\}$, we take the last two numbers $M_{k_{d-1}}, M_{k_d}$ and construct a spline with one knot $\varphi^{1}(t)$ such that $\left\|\left(\varphi^{1}\right)^{(k_i)}\right\| = M_{k_i}$, $i=d-1,d$. Three cases are possible:
\begin{enumerate}
\item If $M_{k_{d-2}} < \left\|\left(\varphi^{1}\right)^{(k_{d-2})}\right\|$ then $M_{{\bf k}}\notin A_d$.
\item If $M_{k_{d-2}} = \left\|\left(\varphi^{1}\right)^{(k_{d-2})}\right\|$ then $M_{{\bf k}} \in A_d \Leftrightarrow M_{{\bf k}} = M_{{\bf k}}(\varphi^{1}).$ 
\item $M_{k_{d-2}} > \left\|\left(\varphi^{1}\right)^{(k_{d-2})}\right\|$.
\end{enumerate}
In the last case there exists a spline $\varphi^2(t)$ with two knots such that $\left\|\left(\varphi^{2}\right)^{(k_i)}\right\| = M_{k_i}$, $i=d-2,d-1,d$. Again three cases are possible:
\begin{enumerate}
\item If $M_{k_{d-3}} < \left\|\left(\varphi^{2}\right)^{(k_{d-3})}\right\|$ then $M_{{\bf k}}\notin A_d$.
\item If $M_{k_{d-3}} = \left\|\left(\varphi^{2}\right)^{(k_{d-3})}\right\|$ then $M_{{\bf k}} \in A_d \Leftrightarrow M_{{\bf k}} = M_{{\bf k}}(\varphi^{2}).$ 
\item $M_{k_{d-3}} > \left\|\left(\varphi^{2}\right)^{(k_{d-3})}\right\|$.
\end{enumerate}
In the last case there exists a spline $\varphi^3(t)$ with three knots such that $\left\|\left(\varphi^{2}\right)^{(k_i)}\right\| = M_{k_i}$, $i=d-3,d-2,d-1,d$. Continuing similarly, we eventually obtain either the needed function that solves Kolmogorov's problem or the fact that $M_{{\bf k}}\notin A_d$.

 {The solution to Kolmogorov's problem in the case when $k_d<r$ is given by the following theorem.}

\begin{theorem}\label{th3'} (Solution to Kolmogorov's problem for $k_d<r$)
{ Let $d\in\NN$, $d\geq 3$, and $0\leq k_1<k_2<...<k_{d}<r$ be
nonnegative integers. Then $ M_{{\bf k}} \in A_{d}(L^{r,r}_{\infty,\infty}({\RR_-}))$ if and only if
$ M_{{\bf {k^2}}} \in A_{d-1}(L^{r,r}_{\infty,\infty}({\RR_-}))$ and
$$M_{k_1}>\lim\limits_{l\to\infty}\left\|\phi_l^{(k_1)}\right\|,$$ where $\phi_l\in \Phi_{r,d-1}$ is such that $\left\|\phi_l^{(k_i)}\right\|=M_{k_i}$, $i = 2,\dots,d$, and $\left\|\phi_l^{(r)}\right\|=l$  $\left(\right.$which exists for all $l\geq \left\|\phi^{(r)}(M_{\bf {k^2}})\right\|$ $\left.\right)$.}
\end{theorem}

{{\bf{Remark.}}}
{It will follow from Lemma~\ref{l1} that $\left\|\phi_l^{(r)}\right\|$ is a decreasing function of $l$.}

\section{Proofs} \label{auxiliary}
In order to prove the main results of this paper, we need several
supporting results, which we present in this section.

{In the case when $k_d=r$} without loss of generality we may assume
that $M_r=1$ and the parameter $l$ in the definition of the spline
is equal to $1$. We also write $\phi(a_1,...,a_s;t)$ instead
of $\phi(a_1,...,a_s,1;t)$.


\begin{lemma}\label{l1}
Let $r,d\in\NN$ and let $0\leq k_1<k_2<...<k_{d}=r$ be integers.
Let a function $x(t)\in L_{\infty,\infty}^{r,r}(\RR)$
 and spline
$\phi\in \Phi_{r,d-2}$ be such that
\begin{equation}\label{rav}
\left\|x^{(k_i)}\right\|=\left\|\phi^{(k_i)}\right\|, \qquad
i=2,...,d.
\end{equation}
 Then
$$
\left\|x^{(k_1)}\right\|\geq\left\|\phi^{(k_1)}\right\|.
$$
Moreover, ``='' is possible only if $x^{(k_1)}\equiv\phi^{(k_1)}$.
\end{lemma}

{\bf Proof of Lemma \ref{l1}.} The lemma will be proved if, assuming
$x^{(k_1)}\ne \phi^{(k_1)}$, we show
$$
\left\|x^{(k_1)}\right\|>\left\|\phi^{(k_1)}\right\|.
$$
Note that the spline $\phi$ has no more than $d-2$ knots. Assume,
contrary to the desired statement, that
\begin{equation}\label{assu}
\left\|x^{(k_1)}\right\|\leq\left\|\phi^{(k_1)}\right\|.
\end{equation}
Denote $\Delta(t):=x(t)-\phi(t)$. To obtain a contradiction, we
will count the number of sign changes of the difference $\Delta(t)$
and its derivatives $\Delta^{(k_1)}(t)$, {$\Delta^{(k_1+1)}(t)$},
..., $\Delta^{(r)}(t)$.

First of all, because of the definition of the spline $\phi$, we
have  $\phi^{(k_1)}(-a_1)=0$. In addition $x^{(k_1)}(-a_1)\geq
0$, and hence we obtain $\Delta^{(k_1)}(-a_1)\geq 0$. Moreover, by assumption
(\ref{assu})
$$
\Delta^{(k_1)}(0)=x^{(k_1)}(0)-\phi^{(k_1)}(0)=\left\|x^{(k_1)}\right\|-
\left\|\phi^{(k_1)}\right\|\leq 0.
$$ Since $x^{(k_1)}\ne
\phi^{(k_1)}$ there exists $t_{k_1+1}^1\in(-a_1,0)$ such that
$\Delta^{(k_1+1)}(t_{k_1+1}^1)< 0$. Besides that,
$\Delta^{(k_1+1)}(-a_1)\geq 0$. Thus, there exists a point
$t_{k_1+2}^1\in(-a_1,0)$ such that
$\Delta^{(k_1+2)}(t_{k_1+2}^1)<0$. Repeating the same arguments, we
obtain that there exist $t_{k_2}^1\in(-a_1,0)$ such that
$\Delta^{(k_2)}(t_{k_2}^1)< 0$. Moreover, $\Delta^{(k_2)}(-a_1)\geq
0$ and $\Delta^{(k_2)}(0)= 0$ by {condition} (\ref{rav}). Therefore,
there exist points $-a_1<t^1_{k_2+1}<t^2_{k_2+1}<0$ such that
$\Delta^{(k_2+1)}(t^1_{k_2+1})<0$ and
$\Delta^{(k_2+1)}(t^2_{k_2+1})>0$. This sign distribution will
remain up to the level $k_3$ where, taking into account {condition}
$\left\|x^{(k_3)}\right\|=\left\|\phi^{(k_3)}\right\|$ and the fact
that $\Delta^{(k_3)}(-a_1)\geq 0$, one can show that there exist
points $-a_1<t^1_{k_3+1}<t^2_{k_3+1}<t^3_{k_3+1}<0$ such that
$\Delta^{(k_3+1)}(t^1_{k_3+1})<0$,
$\Delta^{(k_3+1)}(t^2_{k_3+1})>0$, and
$\Delta^{(k_3+1)}(t^3_{k_3+1})<0$. Proceeding similarly, we obtain
that there exist points
$-a_1<t_{k_{d-1}+1}^1<...<t_{k_{d-1}+1}^{d-1}<0$ such that
$(-1)^i\Delta^{(k_{d-1}+1)}(t_{k_{d-1}+1}^{i})>0$, $i=1,...,d-1$,
and so forth, up to the level $r-1$. Note that ``passing through''
each of the levels $k_1, k_2, ..., k_{d-1}$ we increase the number
of sign changes of the derivatives of the difference function by
one, due to condition (\ref{rav}). At the level $r-1$, there exist
points $-a_1<t_{r-1}^1<...<t_{r-1}^{d-1}<0$ such that
$(-1)^i\Delta^{(r-1)}(t_{r-1}^{i})>0$, $i=1,2,...,d-1$. Moreover,
$\Delta^{(r-1)}(-a_1)\geq 0$. But then on the interval $(-a_1,
t^1_{r-1})$ there exists a set $S_0$ of positive measure such that
$\Delta^{(r)}(t)<0$ for all $t\in S_0$. In addition, for all
$i=1,\dots, d-2$ on the interval $(t^i_{r-1}, t^{i+1}_{r-1})$ there
exists a set $S_i\subseteq (t^i_{r-1}, t^{i+1}_{r-1})$ of positive
measure such that $(-1)^i\Delta^{(r)}(t)<0$ for all $t \in S_i$
(i.e. the function $\Delta^{(r)}(t)$ has no fewer than $d-2$
essential sign changes on $(-a_1, 0)$).

 However, this is impossible, because $(-1)^i\Delta^{(r)}(t)\leq 0$ for almost all $t\in (-a_i, -a_{i+1})$, $i=1,...d-2$ where $a_{d-1}:=0$ and, hence, cannot have more than $d-3$ essential sign changes. We have obtained a
contradiction and, therefore, the lemma is proved. $\square$

{\bf Remark.} {It is easy to see that the statement of the above lemma
remains true if we require
$\left\|x^{(r)}\right\|\leq\left\|\phi^{(r)})\right\|$ instead of
$\left\|x^{(r)}\right\|=\left\|\phi^{(r)}\right\|$.}

From this lemma, as a corollary, we also obtain the following lemma

\begin{lemma}\label{l::2}
Let $r,d\in\NN$ and $0\leq k_1<...<k_{d}=r$ be integers. Let splines
$\phi_1(t)=\phi(a_1,a_2,...,a_s;t)\in {\Phi_{r,d-1}}$ and
$\phi_2(t)=\phi(b_1,b_2,...,b_s;t)\in {\Phi_{r,d-1}}$  be such that
$\left\|\phi_1^{(k_i)}\right\|=\left\|\phi_2^{(k_i)}\right\|$,
$i=1,2,...,d$. Then $\phi_1\equiv \phi_2$.
\end{lemma}

Note that from Lemmas \ref{l1} and \ref{l::2} it follows that the sets of Type~1, Type~2 and Type~3 are pairwise disjoint.

\begin{lemma}\label{l::2'}
Let $d\in\NN$, $d\geq 2$, $x_1>x_2>...>x_d>0$, and
$\alpha_1>\alpha_2>...>\alpha_d$ be given. Then
 $$\begin{vmatrix} x_1^{\alpha_1} & x_2^{\alpha_1}&...&x_d^{\alpha_1}  \\x_1^{\alpha_2} & x_2^{\alpha_2}&...&x_d^{\alpha_2} \\.&.&.&.\\x_1^{\alpha_d} & x_2^{\alpha_d}&...&x_d^{\alpha_d} \end{vmatrix}>0.$$
\end{lemma}
The proof of this statement can be found, for example, in
~\cite{Karlin}. The proof of the next statement is a simple calculus
exercise.


  \begin{lemma}\label{l::4}
 Let $\alpha>\beta>0$ and $M>0$. Then
 $$
\lim\limits_{x\to+\infty}  \left(M+x^\beta\right)^{\frac{\alpha}{\beta}}-x^\alpha=+\infty.$$
\end{lemma}
  As a corollary of it, we obtain

  \begin{lemma}\label{l::sl1}
Let $\alpha>\beta>0$, $M>0$. Let function $y=y(x)$ be such that
$y^\beta-x^\beta\geq M$ for all $x\geq 0$. Then
$y^\alpha-x^\alpha\to+\infty$ as $x\to+\infty$.
  \end{lemma}

The next lemma is very important for the proof of the main
result of this paper.

\begin{lemma}\label{l::th3}
Let $d\in\NN$, $d\geq 3$, and $0\leq k_1<k_2<...<k_{d}=r$ be given.
Let also positive numbers $M_{k_1},...,{M_{k_{d-1}}},M_{k_{d}}$, be
given. In addition, let $M_{{\bf k}^2}\in A_{d-1}^1$ and
$M_{k_1}>\left\|\phi^{(k_1)}(M_{\bf k}^2)\right\|$. Then $M_{{\bf
k}}\in A_{d}^1$.
 \end{lemma}
 {\bf Proof.}
 Without loss of generality we may assume $M_{k_{d}}=1$.

Let us set
$$X:=\left\{a_{d-1}\geq 0\mid \exists
\,a_1>a_2>...>a_{d-2}>a_{d-1}\;\;\right.
$$
$$\left. \hbox{s.t.}\;\;
\left\|\phi^{(k_i)}(a_1,a_2,...,a_{d-1})\right\|=M_{k_i}, \;\;\;
i=2,...,d\right \}.$$ It follows from the assumption $M_{{\bf
k}^2}\in A_{d-2}^1$ that $0\in X$ and, hence, $X\neq \emptyset$.

Let us show that on set $X$ we have well-defined functions $$
a_1(a_{d-1}), \;\; a_2(a_{d-1}),\;\; \dots, \;\;a_{d-2}(a_{d-1}) $$
such that
\begin{equation}\label{st1}
a_1(a_{d-1}) > a_2(a_{d-1} )> \dots > a_{d-2}(a_{d-1})>a_{d-1}
\end{equation}
and
\begin{equation}\label{eq}
\left\|\phi^{(k_j)}(a_1(a_{d-1}),a_2(a_{d-1}),...,a_{d-1}(a_{d-1}),
a_{d-1})\right\|=M_{k_j}, \;\;\; j=2,\dots,d.
\end{equation}
 Existence of numbers $a_1(a_{d-1}), \;a_2(a_{d-1}),\;\dots, \;a_{d-2}(a_{d-1})
$ (for each $a_{d-1}\in X$) with desired properties is obvious (follows from definition of $X$). Let us show that for each $a_{d-1}\in X$ such a set of numbers is unique.

Let us assume that for some $a_{d-1}\in X$ there exist two sets of
knots, namely $a_1>a_2>...>a_{d-2}>a_{d-1}$ and
$b_1>b_2>...>b_{d-2}>a_{d-1}$, such that for $i=2,...,d$
$$\left\|\phi^{(k_i)}(a_1,a_2,...,a_{d-2}, a_{d-1})\right\|=\left\|\phi^{(k_i)}(b_1,b_2,...,b_{d-2}, a_{d-1})\right\|=M_{k_i}.$$ This implies that  $$\left\|\phi^{(k_i)}(a_1,a_2,...,a_{d-2})\right\|=\left\|\phi^{(k_i)}(b_1,b_2,...,b_{d-2})\right\|, \;\;\;\; i=2,...,d,$$ and therefore by Lemma \ref{l::2},  splines are identical. 

Next we show that $X$ is an open subset of $[0, \infty)$ (i.e.
is the intersection of an open set in $\RR$ with $[0, \infty)$) and
that the functions $a_1(a_{d-1})$, $a_2(a_{d-1})$, ..., $a_{d-2}(a_{d-1})$, where $a_{d-1}\in X$, are continuous (and even
smooth).

First, we re-write condition (\ref{eq}) in the following form
\begin{equation} \label{eqq}
F_i(a_1,...,a_{d-2}, a_{d-1})=0, \qquad i=2,...,d-1,
\end{equation}
where functions $F_i$ are defined on the set $a_1>a_2>\dots>a_{d-2}>a_{d-1}$ with $a_1>a_2>\dots>a_{d-2}>0$ (we do not require $a_{d-1}>0$) as follows
$$F_i(a_1,a_2,...,a_{d-1}):=\frac{1}{(r-k_i)!}\sum\limits_{j=1}^{d-1}
(-1)^{j+1} a_j^{r-k_i}-M_{k_i}, \qquad i=2,...,d-1.
$$
By the definition of the set $X$ and the above argument, system (\ref{eqq}) has (for any $a_{d-1}\in X$) a unique solution satisfying (\ref{st1}). In addition, for each set of numbers $a_1>... > a_{d-2} > a_{d-1}$ the Jacobian

$${\rm det}\begin{vmatrix} \frac{\partial F_2}{\partial a_1}&...&\frac{\partial F_2}{\partial a_{d-2}}
 \\.&.&.\\\frac{\partial F_{d-1}}{\partial a_1}&...&\frac{\partial F_{d-1}}{\partial a_{d-2}} \end{vmatrix}
 (a_1, \dots, a_{d-1})\ne 0$$
by Lemma 3.
 Therefore, due to the Implicit Function Theorem (see, for instance, Chapter 8
 in ~\cite{Zorich})
 for any $a_{d-1}\in X$ there exists a neighborhood of the point
 $$\left (a_1(a_{d-1}), \dots, a_{d-2}(a_{d-1}), a_{d-1}\right )$$ of the form
  \begin{equation}\label{111}
\displaystyle \left(\prod_{i=1}^{d-1}(a_i(a_{d-1})-\varepsilon_i,
a_{i}(a_{d-1})+\varepsilon_i)\right)\times
(a_{d-1}-\varepsilon_{d-1}, a_{d-1}+\varepsilon_{d-1}),
 \end{equation}
($\varepsilon_i>0,
a_i(a_{d-1})-\varepsilon_i>a_{i+1}(a_{d-1})+\varepsilon_{i+1},\;\;
i=1,\ldots, d-2, \;\; \varepsilon_{d-1}>0$)
 and continuous (and even smooth) functions
{
$$
b_i\; :\;(a_{d-1}-\varepsilon_{d-1}, a_{d-1}+\varepsilon_{d-1})\to
(a_i(a_{d-1})-\varepsilon_i, a_{i}(a_{d-1})+\varepsilon_i),\;\;
i=1,\ldots , d-2,
$$
such that
$$
b_1(t)>\ldots>b_{d-2}(t),
$$
and for any point $(a_1,\ldots,a_{d-2},t)$ in the neighborhood defined in
(\ref{111})
$$
F(a_1,\ldots,a_{d-2},t)=0
$$
if and only if
$$
a_i=b_i(t),\qquad i=1,\ldots, d-2.
$$
Hence, for any $a_{d-1}\in X$ arbitrary $t\in
(a_{d-1}-\varepsilon_{d-1}, a_{d-1}+\varepsilon_{d-1})$ will belong
to $X$. This implies that $X$ is an open set and functions
$a_i(a_{d-1}),\; i=1,\ldots, d-1,$ are continuous on $X$.}

Let us take the maximal connected component of the set $X$ that contains $0$.
We next show that it is bounded.

Assuming the contrary, we consider the sequence $\{
a^n_{d-1}\}_{n=1}^\infty$ of the points of this component such that
$a^n_{d-1}\to \infty$ as $n\to \infty$. Let us also consider the
corresponding sequences $\{ a_i(a^n_{d-1})\}_{n=1}^\infty$,
$i=1,\dots,d-2$. It is clear that $\{
a_i(a^n_{d-1})\}_{n=1}^\infty\to \infty$ as $n\to \infty$,
$i=1,\ldots,d-2,$ as well.

Let us consider {the equalities}
\begin{equation}\label{st}
\sum\limits_{j=1}^{d-1} (-1)^{j+1}
(a_j(a^n_{d-1}))^{r-k_i}=(r-k_i)!M_{k_i}, \qquad i=2,...,d-1,
\end{equation}
{(here $a_{d-1}(a_{d-1}):=a_{d-1}$).}

Switching, if needed, to a subsequence, we see that for $i=d-1$
there exists an index $j$ and constant $c>0$ such that { for all
$n\in\mathbb{N}$}
$$
(a_{2j-1}(a^n_{d-1}))^{r-k_{d-1}}-(a_{2j}(a^n_{d-1}))^{r-k_{d-1}}{>c(r-k_{d-1})!M_{k_{d-1}}}.
$$
Then, by Lemma 5, for all $i<d-1$ we have
$$
(a_{2j-1}(a^n_{d-1}))^{r-k_{i}}-(a_{2j}(a^n_{d-1}))^{r-k_{i}}\to
\infty, \qquad n\to \infty,
$$
and, therefore, for all $i<d-1$
$$\sum\limits_{j=1}^{d-1}
(-1)^{j+1} (a_j(a^n_{d-1}))^{r-k_i}\to\infty, \qquad n\to\infty,
$$
which contradicts the fact that (\ref{st}) holds for all $i=2,\dots,d-1$.

Thus, the maximal connected component of the set $X$, that contains zero, is indeed bounded. Let us denote its right end-point by $a_{d-1}^*$. Clearly, $a_{d-1}^*\notin X$.

Next we show that
\begin{equation}\label{3}
\lim\limits_{a_{d-1}\to a_{d-1}^*}
\left\|\phi^{(k_1)}(a_1(a_{d-1}),...,a_{d-2}(a_{d-1}),a_{d-1})\right\|=\infty.
\end{equation}

To do so, we first show that the set $\left\{a_2(a_{d-1})\mid
a_{d-1}\in[0,a_{d-1}^*)\right\}$ is unbounded.
 Assume to the contrary that it is bounded. It implies that the set $\left\{a_1(a_{d-1})\mid a_{d-1}\in[0,a_{d-1}^*)\right\}$ is bounded as well.

In addition, the following sets are also bounded
$$
\left \{ a_i(a_{d-1}) \; : \; \; a_{d-1}\in [0, a_{d-1}^*)\right \}, \qquad i=3, \dots, d-2.
$$
Let us take an arbitrary sequence $a_{d-1}^n\to a_{d-1}^*$ as $n\to \infty$. Switching to a subsequence, if necessary, we may assume that all subsequences $\{a_i(a_{d-1}^n) \}_{n=1}^{\infty}$ have limits. Set
$$
a_i^*:= \displaystyle \lim_{n\to \infty} a_i(a_{d-1}^n), \qquad i=1, \dots, d-2.
$$
Clearly, {$a_1^*\geq a_2^*\geq ...\geq a_{d-2}^*\geq a_{d-1}^*$}.

Since for all $n$ the vector
$(a_1(a_{d-1}^n),...,a_{d-2}(a_{d-1}^n),a_{d-1}^n)$ is the solution
of the system
\begin{equation}\label{1}
F_i(a_1,a_2,...,a_{d-1})=0,\qquad i=2,...,d-1,
\end{equation}
then the vector $(a_1^*,...,a_{d-2}^*,a_{d-1}^*)$ is also a solution of
the system {\eqref{1}}. In addition, since $a_{d-1}^*\notin X$, there
exists $i\in\left\{1,2,...,d-2\right\}$ such that $a_i^*=a_{i+1}^*$
(otherwise we would have $a_{d-1}^*\in X$). This eliminates at least two knots of the spline. Then there exists a spline
$\phi\in\Phi_{r, d-3}$, such that
$$
\left\|\phi^{(k_s)}\right\|=\left\|\phi^{(k_s)}(a_1(0),...,a_{d-2}(0),0)\right\|,
\qquad s=2,...,d.
$$
Indeed, it is sufficient to take the spline with different knots out
of the set $a_1^*,...,a_{d-2}^*,a_{d-1}^*$. However, it contradicts
Lemma \ref{l1}, which implies that when the corresponding norms are
equal for the splines $\phi \in \Phi_{r,d-2}$ and $\phi \in \Phi_{r,
d-3}$ we obtain
$$\left\|\phi^{(k_2)}\right\|<\left\|\phi^{(k_2)}(a_1(0),...,a_{d-2}(0),0)\right\|.$$
Hence, we proved that the set $\left\{a_2(a_{d-1})\mid
a_{d-1}\in[0,a_{d-1}^*)\right\}$ is unbounded.

If for some $\varepsilon>0$ and some $s\in\left\{2,3,...,{d-1}\right\}$
we have $a_1^{r-k_s}(a_{d-1})-a_2^{r-k_s}(a_{d-1})>\varepsilon$ in some
neighborhood $(a_{d-1}^*-\delta,a_{d-1}^*)$, then by Lemma
\ref{l::sl1} we have $a_1^{r-k_1}(a_{d-1})-a_2^{r-k_1}(a_{d-1})\to\infty$ as
$a_{d-1}\to a_{d-1}^*$ and, hence, (\ref{3}) is proved.

Let us assume that for any $s=2,\dots, d-1$ we have $a_1^{r-k_s}(a_{d-1})-a_2^{r-k_s}(a_{d-1})\to {0}$ as
$a_{d-1}\to a_{d-1}^*$. Next we show that the set
$$
\left\{ a_4(a_{d-1})\;: \;\; a_{d-1}\in [0, a_{d-1}^*)\right \}
$$
is unbounded. Assuming the contrary, we choose a sequence $a_{d-1}^n\to a_{d-1}^*$ as $n\to \infty$ and, switching if necessary to subsequences, we can assume that all sequences $\{a_i(a_{d-1}^n) \}_{n=1}^{\infty}$ have limits $a_3^*\geq a_4^*\geq \dots \geq a_{d-2}^*\geq a_{d-1}^*$.

In addition, since for any $s=2, \dots, d-1$ we have
$$a_1^{r-k_s}(a_{d-1})-a_2^{r-k_s}(a_{d-1})\to 0, \qquad a_{d-1}\to a_{d-1}^*,
$$
and numbers $a_1^n(a_{d-1}), \dots, a_{d-2}^n(a_{d-1}), a_{d-1}^n$ solve the system of equations (\ref{st}). We take the limit as $n\to \infty$ to obtain
$$
\displaystyle \sum_{j=3}^{d-1}(-1)^{j+1}(a_j^*)^{r-k_i}={(r-k_i)!M_{k_i}}, \qquad i=2, \dots, d-1.
$$
Thus, there exists a spline $\phi\in \Phi_{r,{ d-3}}$ such that
$$
\|\phi^{(k_i)}\|=\|\phi^{(k_i)}(a_1(0), \dots, a_{d-2}(0), 0)\|,
\qquad i=2,\dots, d-1.
$$
However, this, as in the case of proving unboundedness of the set $$\left\{ a_2(a_{d-1})\;: \;\; a_{d-1}\in [0, a_{d-1}^*)\right \},$$ contradicts Lemma \ref{l1}. Therefore, the set $\left\{ a_4(a_{d-1})\;: \;\; a_{d-1}\in [0, a_{d-1}^*)\right \}$ is unbounded.

Next, if for some $\varepsilon >0$ and some $s\in \{2,3,\dots, d-1\}$ we have
$$
a_3^{r-k_s}(a_{d-1})-a_4^{r-k_s}(a_{d-1})>\varepsilon
$$
in some neighborhood $(a^*_{d-1}-\delta, a^*_{d-1})$, then, as above,
we establish that the function $$\| \phi^{(k_1)}(a_1(a_{d-1}),
\dots, a_{d-2}(a_{d-1}), a_{d-1}) \|$$ is unbounded and this
concludes the proof of {(\ref{3}).}

If for {all} $s\in \{2,3,\dots, d-1\}$ we have
$$
a_3^{r-k_s}(a_{d-1})-a_4^{r-k_s}(a_{d-1})\to 0, \qquad a_{d-1}\to a_{d-1}^*,
$$
then similarly to the above, we establish that the set $$\left\{ a_6(a_{d-1})\;: \;\; a_{d-1}\in [0, a_{d-1}^*)\right \}$$ is unbounded.


Continuing in a similar way, we either find such
$s\in\left\{2,3,...,d-1\right\}$ and $i\in
\left\{1,2,...,\left[\frac{d-1}{2}\right]\right\}$, that
$$a_{2i-1}^{r-k_s}-a_{2i}^{r-k_s}>\varepsilon>0\qquad \hbox {in some
neighborhood} \;\;\; (a_{d-1}^*-\delta,a_{d-1}^*)$$ (which by Lemma
\ref{l::sl1} implies (\ref{3})), or we obtain that $d-1$ is odd and
$$(a_{d-1}^*)^{r-k_s}=M_k(r-k_s)!\qquad  \hbox{for all} \;\;\;
s=2,...,d-1.$$ However, in the latter case the spline
$\phi(a_{d-1}^*;t)$ satisfies
$$
\left\|\phi^{(k_s)}(a_{d-1}^*)\right\|=\left\|\phi^{(k_s)}(a_1(0),...,a_{d-2}(0),0)\right\|,
\qquad s=2,...,{d},
$$ which contradicts Lemma \ref{l1}.

Therefore, we have proved (\ref{3}). Hence, due to the continuity
of
$$\left\|\phi^{(k_1)}(a_1(a_{d-1}),...,a_{d-2}(a_{d-1}),a_{d-1})\right\|,$$
there exists  $a_{d-1}>0$ such that
$\left\|\phi^{(k_1)}(a_1(a_{d-1}),...,a_{d-2}(a_{d-1}),a_{d-1})\right\|=M_{k_1}$
which concludes the proof of the lemma. $\square$

   \begin{lemma}\label{l::th2}
Let $d\in\NN$, $d\geq 3$, integers $0\leq k_1<k_2<...<k_{d}=r$,
and positive real numbers $M_{k_1},...,M_{k_{d-1}},M_{k_{d}}$ be
given. Then
$$
\left \{ M_{{\bf k}}\in A_d(L^{r,r}_{\infty,\infty}({\RR_-}))\right \} \;\; \Longleftrightarrow \;\; \left \{ M_{{\bf k}} \in A_d^1\cup A_d^2\cup A_d^3\right \}.
$$
 \end{lemma}
 {\bf Proof.}
 The sufficiency part is obvious.

 Let us prove the necessity. Again, without loss of generality, we may assume $M_d=1$.
 We proceed by induction on $d$.
For $d=3$ the statement follows from the results of Olovyanishnikov
($k_1=0$) and V. Babenko and Y. Babenko ($k_1>0$). This is the basis of induction.
Let $m\geq 4$. Assume that the statement is true for $3\leq d \leq m-1$. We prove that from here it follows that the statement is also true for $d=m$.

 Let
$M_{{\bf k}}\in A_m(L^{r,r}_{\infty,\infty}({\RR_-}))$ and let $x\in L_{\infty,\infty}^{r,r}(\RR_-)$ be such
that
$$\left\|x^{(k_i)}\right\|=M_{k_i}, \qquad i=1,...,m.
$$
 Then the
set $M_{{\bf k}^2}\in A_{m-1}(L^{r,r}_{\infty,\infty}({\RR_-}))$
and, hence, by induction assumption there exists a spline
$\phi(t)=\phi(a_1,...,a_s;t)\in\Phi_{r, m-2}$, such that
$$
\left\|x^{(k_i)}\right\|=\left\|\phi^{(k_i)}\right\|, \qquad
i=2,...,m.
$$
 By Lemma \ref{l1} we have
$\left\|x^{(k_1)}\right\|\geq\left\|\phi^{(k_1)}\right\|$. If
$\left\|x^{(k_1)}\right\|=\left\|\phi^{(k_1)}\right\|$, then
$M_{{\bf k}}\in A^2_{m}(L^{r,r}_{\infty,\infty}({\RR_-}))$. Let
 \begin{equation}\label{16}
 \left\|x^{(k_1)}\right\|>\left\|\phi^{(k_1)}\right\|.
 \end{equation}
If $k_1=0$, then there exists $C>0$ such that
$\left\|(\phi+C)^{(k_i)}\right\|=M_{k_i}$, $i=1,2,...,m$, and hence
the set $M_{{\bf k}}\in A_d^3{\bigcup A_d^1}$.

Let now $ k_1>0$.
 Then by Lemma \ref{l1} we have $s=m-2$. Indeed, otherwise (if $s<m-2$) from $$\left\|x^{(k_{i})}\right\|=\left\|\phi^{(k_{i})}\right\|, \qquad i=3,...,m$$ by Lemma \ref{l1} (with $d=m-1$) it would follow that
\begin{equation}\label{18}
 \left\|x^{(k_2)}\right\|\geq\left\|\phi^{(k_2)}\right\|,
\end{equation}
with ``='' possible only when $x^{(k_2)}\equiv \phi^{(k_2)}$. Since
$x,\phi\in L_{\infty,\infty}^{r,r}(\RR_-)$, and due to the fact that
$k_1>0$, from  $x^{(k_2)}\equiv \phi^{(k_2)}$ it would follow that
$x^{(k_1)}\equiv \phi^{(k_1)}$, which would contradict \eqref{16}.

Thus, $M_{{\bf k}^2}\in A^1_{m-1}$. Taking into account \eqref{16},
we obtain $$M_{k_1}>\left\|\phi^{(k_1)}\right\|.$$ Therefore, by
Lemma \ref{l::th3}, the set $M_{{\bf k}}\in A_{m}^1$ which concludes
the step of induction. $\square$

\begin{lemma}\label{l::8}
{Let $d\in\NN$, $d\geq 3$, integers $0\leq k_1<k_2<...<k_{d}<r$, and
a function $x\in L_{\infty,\infty}^{r,r}(\RR_-)$ be given. Then
there exists a spline $\phi\in\Phi_{r,d-1}$ such that }
\begin{equation}\label{l8.1}
{\left\|\phi^{(k_i)}\right\|=\left\|x^{(k_i)}\right\|,\qquad
i=1,\dots,d.}
\end{equation}
{Moreover, if $\phi^{(k_1)}\neq x^{(k_1)}$ then  }
\begin{equation}\label{l8.2}
{\left\|\phi^{(r)}\right\|<\left\|x^{(r)}\right\|.}
\end{equation}
\end{lemma}
{\bf Proof.}
{Consider the following extremal problem
$$\left\|\phi^{(r)}\right\|\to {\rm min},$$}
over all splines
$\phi\in\Phi_{r,d}$ which satisfy conditions \eqref{l8.1}. First of
all note that due to Lemma~\ref{l::th2} the set of splines which
satisfy conditions \eqref{l8.1} is nonempty. It is easy to see
that the minimum in the above extremal problem is attained on some spline
$\phi^*\in\Phi_{r,d}$. Next we will show that $\phi^*\in\Phi_{r,d-1}$. Assume to
the contrary, that the spline $\phi^*$ has exactly $d$ knots
$-a_1^*,\dots,-a_{d}^*$. Consider the following Lagrange function
{$$L:=l+\sum_{i=1}^d\lambda_i\left(\left\|\phi^{(k_i)}\right\|-M_{k_i}\right)=l+\sum_{i=1}^d\lambda_i\left[\frac{l}{(r-k_i)!}\sum_{j=1}^d(-1)^{j+1}a_j^{r-k_i}-M_{k_i}\right],$$}
{where $l=\left\|\phi^{(r)}\right\|$,
$M_{k_i}=\left\|x^{(k_i)}\right\|$ and $-a_1,\dots,-a_d$ are the
knots of the spline $\phi$. Due to the necessary condition for
optimality there exist numbers $\lambda_1,
\lambda_2,\dots,\lambda_d$ such that}
\begin{equation}\label{l8.3}
{\sum_{i=1}^d\lambda_i\left[\frac{l^*}{(r-k_i-1)!}(-1)^{j+1}(a_j^*)^{r-k_i-1}\right]=0,\qquad j=1,\dots,d,}
\end{equation}
{and}
\begin{equation}\label{l8.4}
{1+\sum_{i=1}^d\lambda_i\left[\frac{1}{(r-k_i)!}\sum\limits_{j=1}^{d}(-1)^{j+1}(a_j^*)^{r-k_i}\right]=0,}
\end{equation}
{where $l^*=\left\|{\phi^*}^{(r)}\right\|$. System \eqref{l8.3} is a system of linear equations with respect to $\lambda_i$, $i=1,2,\dots,d$. The main determinant of system \eqref{l8.3} is nonzero due to Lemma~\ref{l::2'}, so it has only trivial solution. But it contradicts condition~\eqref{l8.4}. So $\phi^*\in\Phi_{r,d-1}$. Due to Lemma~\ref{l1}, inequality~\eqref{l8.2} holds.  Lemma is proved. $\square$}

\begin{lemma}\label{l::9}
{Let $d\in\NN$, $d\geq 2$, integers $0\leq
k_1<k_2<...<k_{d}<r$ and a function $x\in L_{\infty,\infty}^{r,r}$
be given. Then for all $l>\left\|x^{(r)}\right\|$ there exists a
spline $\phi\in\Phi_{r,d}\setminus \Phi_{r,d-1}$ such that
$\left\|\phi^{(k_i)}\right\|=\left\|x^{(k_i)}\right\|,\,i=1,\dots,d$,
and $\left\|\phi^{(r)}\right\|=l$.}
\end{lemma}
{We will prove the Lemma by induction on $d$. The
case when $d=2$ follows from Olovyanishnikov inequality~\eqref{olov}
and the work of V.Babenko and Y.Babenko~\cite{us}. Assume that the statement of the Lemma holds for some  $d \in\NN$, $d\geq 2$. We will prove that it also holds
for $d+1$. Let integers $0\leq k_1<k_2<...<k_{d+1}<r$,  a
function $x\in L_{\infty,\infty}^{r,r}$ and
$l>\left\|x^{(r)}\right\|$ be given. By induction hypothesis, there
exists a spline $\phi\in \Phi_{r,d}\setminus \Phi_{r,d-1}$ such that
$\left\|\phi^{(k_i)}\right\|=\left\|x^{(k_i)}\right\|,\,i=2,3,\dots,d
+1$ and $\left\|\phi^{(r)}\right\|=l$. Due to  Lemma~\ref{l1} we
have $\left\|\phi^{(k_1)}\right\|<\left\|x^{(k_1)}\right\|$.
Applying Lemma~\ref{l::th3}, we obtain that the statement of
Lemma~\ref{l::9} holds true.} $\square$

Finally, let us turn to the proofs of the main results. First of all observe that {existence of the spline in} Theorem \ref{th1} follows from Lemma \ref{l::th2} { and Lemma~\ref{l::8}. The proof of the extremal properties for the case $k_d=r$ is similar to the proof of Lemma~\ref{l1}, and the proof of the extremal property for the case $k_d<r$ is contained in Lemma~\ref{l::8}.} As we already mentioned earlier, Theorem~\ref{th1} can be viewed as a generalization of Olovyanishnikov's inequality to the case of $d> 3$ consecutive derivatives of an arbitrary function $x\in L_{\infty,\infty}^{r,r}(\RR_-)$.

Theorem \ref{T3} immediately follows from Lemmas \ref{l::th2} and \ref{l::8}. { Theorem \ref{th3} follows from combining Theorem~\ref{th1} and Lemmas \ref{l1}, \ref{l::th3}, and \ref{l::th2}. Finally, Theorem~\ref{th3'} follows from Lemmas~\ref{l1}~and~\ref{l::9}.}

\end{document}